\documentclass[12pt]{article}

\usepackage{latexsym,amssymb,amsmath}
\usepackage{amsthm}
\renewcommand{\qed}{\hfill\rule[-0pt]{5pt}{15pt}}

\font\defi=cmr7
\newcommand{\adef}{\,\,\textrm{\raisebox{1.5pt}{\defi :}=}\,}

\newcommand{\aeq}{\,\,\textrm{\phantom{\raisebox{1.5pt}{\defi :}}=}\,}

\usepackage{amsopn}
\DeclareMathOperator{\re}{Re}

\newtheorem{lemma}{Lemma}
\newtheorem{conjecture}{Conjecture}

\newcommand{\nn}{\mathbb N}

\newcommand{\rr}{\mathbb R}
\newcommand{\cc}{\mathbb C}
\newcommand{\kk}{\mathbb K}
\newcommand{\mm}{\mathbb M}
\newcommand{\modulo}[1]{\left\vert{#1}\right\vert}
\newcommand{\norma}[1]{\left\Vert{#1}\right\Vert}
\newcommand{\pent}[1]{\left[{#1}\right]}
\newcommand{\pfrac}[1]{{#1}-\pent{{#1}}}
\newcommand{\frho}[1]{\rho\left({#1}\right)}
\newcommand{\bin}[2]{\frho{\frac{x}{{#1}}}\;-\;\frac{{#2}}{{#1}}\,\frho{\frac{x}{{#2}}}}

\begin{document}


\begin{center}

{\bf\Large A Sequence of Beurling Functions Related to the Natural
Approximation $B_{n}$ Defined by an Iterative Construction Generating
Square-Free Numbers $k_{i}$ and the Value of the M\"obius Function
$\mu(k_{i})$}

\vspace{5mm}

F. Auil

\end{center}

\vspace{0.2cm}


\begin{quote}
{\small
{\bf Abstract.} For a function 
$F_{n}(x)=\sum_{k=1}^{n}a_{k}\,\rho\left(\frac{\theta_{k}}{x}\right)$, 
where $\rho(x)\adef\pfrac{x}$, with $a_{k}\in\cc$ and 
$0<\theta_{k}\leqslant 1$ satisfying
$\sum_{k=1}^{n}a_{k}\,\theta_{k}=0$, is known to have 
$$
\modulo{\frac{1}{s}\left(1-\zeta(s)\,\sum_{k=1}^{n}a_{k}\,\theta_{k}^{s}\right)}
=\modulo{\int_{0}^{1} (F_{n}(x)+1)\,x^{s-1}\,dx}
\leqslant\norma{F_{n}(x)+1}\,\norma{x^{s-1}};
$$
where the first relation follows from straightforward computation and
the second using Schwarz inequality in $L^{2}([0,1],dx)$. Therefore,
if the first norm in the right hand side above would be arbitrarily
small for a suitable choice of $n$, $a_{k}$'s and $\theta_{k}$'s the
function $\zeta(s)$ would have zeros for $\re\,s>1/2$. This is the
Beurling approach to Riemann Hypothesis. Several \emph{approximating
sequences} were proposed, between them 
$$
B_{n}(x)\adef\sum_{k=1}^{n} \mu(k)\,\frho{\frac{1/k}{x}} \;
- \;
n\left(\sum_{k=1}^{n}\frac{\mu(k)}{k}\right)\frho{\frac{1/n}{x}}. 
$$
In the present work we construct \emph{iteratively} a sequence of numbers
$\{k_{n}\}$ and approximating functions $\{\widetilde{B}_{n}\}$ converging
\emph{pointwise} to $-1$ in $[0,1]$. 
We prove results which seems to suggest the relation
$\widetilde{B}_{n}={B}_{n}$ and we conjecture that a
sufficient condition for this is that the set $\{k_{n}\}_{n\in\nn}$ be
equal to the set of square-free numbers,
$\mm\adef\{m\in\nn\,:\,\mu(m)\neq0\}$. 
Numerical evidence seems to support both conjectures. 

Anyway, we think that these sequences are interesting by itself because
our construction not only generates square-free (hence prime)
numbers $k_{i}$, but also the value of the M\"obius function
$\mu(k_{i})$. Our definition is independent of $\mm$ and $\mu$, with the
$k_{i}$'s arising as discontinuity points of the
$\widetilde{B}_{i}$'s. 

As for the case of $B_{n}$, we prove that sequence
$\widetilde{B}_{n}$ is not convergent to $-1$ in
$L^{2}([0,1],dx)$. Consequently, we focus our analysis not on 
$L^{2}$ norm at right hand side in first expression above but on the
integral at the middle term. This procedure seems to be useful to
elucidate the lack of $L^{2}$ convergence for \emph{step} Beurling
functions.
} 
\end{quote}


\section{Introduction: Beurling Functions and Approximating Sequences}

Denote as $\pent{x}$ the integer part of $x$, i.e. the greatest
integer less than or equal to $x$ and define the fractional part 
function by $\rho(x)=\pfrac{x}$. 
Given $n\in\nn$ and two families of parameters
$\{a_{k}\}_{k=1}^{n}\subset\cc$ and $\{\theta_{k}\}_{k=1}^{n}\subset
(0,1]$, we define a \emph{Beurling function} as a function
$F=F_{n}$ (the sub-index $n$ included in the notation for
convenience) of the form
\begin{equation}
F_{n}(x)\adef\sum_{k=1}^{n} a_{k}\,\frho{\frac{\theta_{k}}{x}}, 
\label{eq-1}
\end{equation} 

\noindent
For a Beurling function $F_{n}$, an elementary computation shows 
\begin{equation}
\int_{0}^{1} (F_{n}(x)+1) \, x^{s-1} \, dx = \frac{\sum_{k=1}^{n} a_{k}\,\theta_{k}}{s-1}
+ \frac{1}{s}\left(1-\zeta(s) \, \sum_{k=1}^{n}
  a_{k}\,\theta_{k}^{s}\right);\quad\re s>0.  
\label{eq-2}
\end{equation}

\noindent
See, for instance, \cite[p. 253]{don} for a proof. 
It is useful, but not always necessary, assume that the parameters
defining the function $F_{n}$ satisfy the additional condition 
\begin{equation}
\sum_{k=1}^{n} a_{k}\,\theta_{k}=0. 
\label{eq-3}
\end{equation}

\noindent
In this case, first term at right hand side of
(\ref{eq-2}) vanishes, simplifying the expression. 
Identity (\ref{eq-2}) is the starting point of a theorem by
Beurling; see \cite[p. 252]{don} for a proof and further references.
Here we just remark that an elementary computation, using Schwarz
inequality for the integral at left hand side of (\ref{eq-2}), allows
to show that a sufficient condition for Riemann Hypothesis (RH) is that
$\norma{F(x)+1}$ be done arbitrarily small for a suitable choice of
$n$, $a_{k}$'s and $\theta_{k}$'s, where $\norma{.}$ denotes the norm
in $L^{2}([0,1],dx)$. We will refer to this last condition as to the
\emph{Beurling criterion} (BC) for RH. 
It was proved in \cite{bae2} that BC remains unchanged if we restrict
the parameters $\theta_{k}$ to be reciprocal of natural numbers,
i.e. $\theta_{k}=1/b_{k}$, with $b_{k}\in\nn$.\\

Several \emph{approximating} functions (\ref{eq-1}) were proposed. 
From (\ref{eq-2}) and (\ref{eq-3}), we have that under BC the
``partial sum''
\begin{equation}
\sum_{k=1}^{n} a_{k}\,\theta_{k}^{s}.
\end{equation}
is an approximation to the inverse Riemann Zeta Function 
$1/\zeta(s)$, which is known to have an expression as a Dirichlet
series  
\begin{equation}
\frac{1}{\zeta(s)}\;=\;\sum_{k=1}^{n} \frac{\mu(k)}{k^{s}}, 
\end{equation}
convergent for $\re s>1$. Therefore, a (naive) first choice for an
approximating function would be
\begin{equation}
S_{n}(x)\adef\sum_{k=1}^{n}\mu(k)\,\frho{\frac{1/k}{x}}. 
\label{eq-6}
\end{equation}
But this function does not matches the condition (\ref{eq-3}). We can
handle this without subtlety, just taking out the difference,
given by $g(n)$, where
\begin{equation}
g(t)\adef \sum_{\nn\ni k\leqslant t}\frac{\mu(k)}{k}. 
\end{equation}
Therefore a second choice would be
\begin{align}
B_{n}(x) & \adef\sum_{k=1}^{n} \mu(k)\,\frho{\frac{1/k}{x}} \;
- \;n\,g(n)\,\frho{\frac{1/n}{x}}.\\
 & \aeq\sum_{k=1}^{n-1} \mu(k)\,\frho{\frac{1/k}{x}} \;
- \;n\,g(n-1)\,\frho{\frac{1/n}{x}}. 
\label{eq-7}
\end{align}
There are also variants on the same theme, as
\begin{equation}
V_{n}(x)\adef\sum_{k=1}^{n} \mu(k)\,\frho{\frac{1/k}{x}} \;
- \; g(n)\,\frho{\frac{1}{x}}.
\label{eq-8}
\end{equation}
Sequences (\ref{eq-6}), (\ref{eq-7}) and (\ref{eq-8}) are known to
be \emph{not} convergent to $-1$ in $L^{2}([0,1],dx)$, as proved in
\cite{bae}.


\section{From Integrals to Series}

The integral in (\ref{eq-2}) can be expressed alternatively as a series,
by a ``change of variable'' under suitable hypothesis. A
Beurling function constant between the reciprocal of the natural
numbers will be called a \emph{step Beurling function}, i.e. such a
function takes (non-necessarily different) constant values in each of
the intervals $(\frac{1}{k+1},\frac{1}{k}]$, for all $k\in\nn$. Note
that Beurling functions are left-continuous, because $\rho(x)$ is
right-continuous.

\begin{lemma} 
Let $F_{n}$ be a step Beurling function, such that
$F_{n}(x)=-1$ if $x\in(\frac{1}{m},1]$, where $m\in\nn$. (If $m=1$
this is an empty condition, therefore in this case there is not
additional condition at all). 
Define $f_{n}(k)\adef F_{n}(1/k)$, for $k\in\nn$. Then, 
\begin{equation}
s\int_{0}^{1} (F_{n}(x)+1) \, x^{s-1} \, dx \;
= \; \sum_{k=m}^{\infty}f_{n}(k)
\left[\frac{1}{k^{s}}-\frac{1}{(k+1)^{s}}\right] \; 
+ \;\frac{1}{m^{s}}.
\label{eq-9}
\end{equation}

\end{lemma}

\paragraph{\it Proof:}
\begin{multline*}
\int_{0}^{1} (F_{n}(x)+1) \, x^{s-1} \, dx \;
= \; \int_{0}^{\frac{1}{m}} (F_{n}(x)+1) \, x^{s-1} \, dx \; 
= \; \sum_{k=m}^{\infty}\int_{\frac{1}{k+1}}^{\frac{1}{k}} (F_{n}(x)+1) \, x^{s-1} \,
dx\\
= \sum_{k=m}^{\infty}(F_{n}(1/k)+1)\int_{\frac{1}{k+1}}^{\frac{1}{k}} x^{s-1} \, dx \; 
= \;
\sum_{k=m}^{\infty}(F_{n}(1/k)+1)\left.\frac{x^{s}}{s}\right\vert_{\frac{1}{k+1}}^{\frac{1}{k}}\\
=\frac{1}{s}\sum_{k=m}^{\infty}(F_{n}(1/k)+1)
\left[\frac{1}{k^{s}}-\frac{1}{(k+1)^{s}}\right] \; 
= \; \frac{1}{s}\sum_{k=m}^{\infty}F_{n}(1/k)
\left[\frac{1}{k^{s}}-\frac{1}{(k+1)^{s}}\right] \; + \;
\frac{1}{s}\frac{1}{m^{s}}.\qed
\end{multline*}\\

\paragraph{Remarks.} 
{\bf (1)} Observe that a step Beurling function $F_{n}(x)$ is
completely determined for $x\in[0,1]$ by the
arithmetic function $f_{n}(k)$, $k\in\nn$. 
This arithmetic function can be extended, just by
defining $f_{n}(x)\adef F_{n}(1/x)$, for $x\in\rr$, becoming
right-continuous. 

\medskip
\noindent
{\bf (2)} Just for reference, we define an 
\emph{arithmetic Beurling function} as a right-continuous 
function $f$ on $[1,+\infty)$ constant between the natural numbers,
such that $f(1/x)$ is a Beurling function. Thus, there exists a
correspondence between step and arithmetic Beurling functions and the
integrals involving the former correspond to series involving the
later as expressed in (\ref{eq-9}). 

\medskip
\noindent
{\bf (3)} We can think in relation (\ref{eq-9}) as a ``change of
variables'' in the integral, turning the integration domain from
$[0,1]$ to $[1,+\infty)$. This can be visualized also like to put a zoom
on the original integration domain, reflecting the interval $[0,1]$
and then stretching it to fit on $[1,+\infty)$.


\section{An Arithmetic Beurling Function Iteratively Defined}
\label{rep}


\subsection{Beurling Binomials}

One of the simplest arithmetic Beurling functions
matching (\ref{eq-3}) are given by
\begin{equation}
\beta_{a,b}(x)\adef\bin{a}{b}, 
\end{equation}
when $a,b\in\nn$. These ``binomials'' will be the basic blocks in our
construction, thus we summarize some of its elementary properties in
the following result. 

\begin{lemma} 
Consider $a,b\in\rr$, with $0<a<b$. Then, 

\begin{itemize}

\item[\bf{a.}] $\displaystyle{\frho{\frac{x}{a}}}$ and 
$\displaystyle{\frho{\frac{x}{b}}}$ are right-continuous, and
linearly independent functions. 

\item[\bf{b.}] 
$\beta_{a,b}(x)=0$, when $0\leqslant x<a$. 

\item[\bf{c.}] Let $k\in\nn$ be such that $(k-1)a<b\leqslant ka$. Then,  
$$
\beta_{a,b}(x)=\begin{cases}-j & \text{\rm if $ja\leqslant x<(j+1)a$, for
$j=1,\dots,(k-2)$;}\\-(k-1) & \text{\rm if $(k-1)a\leqslant x<b$.}\end{cases}
$$

\item[\bf{d.}] Assume $a,b\in\nn$. Then, 
$\beta_{a,b}(x)$ is constant when $k\leqslant x<k+1$, for all
$k\in\nn$. 

\end{itemize}
\label{lemma-1}
\end{lemma}


\subsection{The Aproximating Sequence $\widetilde{B}_{i}$}

We will define a sequence of numbers $\{k_{i}\}$ and functions
$\{b_{i}\}$ iteratively as follows. We start with the definition
\begin{equation}
\begin{split}
k_{1} & \adef 1;\\
k_{2} & \adef 2;\\
b_{2}(x) & \adef\frho{\frac{x}{k_{1}}}\; 
- \; \frac{k_{2}}{k_{1}}\,\frho{\frac{x}{k_{2}}}. 
\label{eq-20}
\end{split}
\end{equation}
\noindent
And for $i\geqslant 2$ define $k_{i+1}\adef k_{i}+j$, where $j$ is
the less integer such that $b_{i}(k_{i}+j)\neq b_{i}(k_{i})$, and 
\begin{equation}
b_{i+1}(x)\adef b_{i}(x) \; 
+ \; \left(1+b_{i}(k_{i})\right)\left[\frho{\frac{x}{k_{i}}}\; 
- \; \frac{k_{i+1}}{k_{i}}\,\frho{\frac{x}{k_{i+1}}}\right]. 
\label{eq-21}
\end{equation}\\

\noindent
Observe that each $b_{i}$ is a linear combination of $\beta_{p,q}$. 
Other elementary properties are given in the next result, which is a
direct consequence of Lemma \ref{lemma-1}. 

\begin{lemma} For any $i\in\nn$ we have
\begin{itemize}

\item[{\bf a.}] $b_{i}$ is an arithmetic Beurling function, i.e. a
right-continuous function constant between the natural numbers. 

\item[{\bf b.}] $b_{i+1}(k_{i})=-1$.

\item[{\bf c.}] Assume $k_{i+1}\leqslant 2k_{i}$ for $i\geqslant 2$. Then,
$b_{i}(x)=-1$ for all $x\in[1,k_{i})$. In particular, the sequence 
$\{b_{i}\}_{i\in\nn}$ converges pointwise to $-1$ in
$[1,+\infty)$. 

\end{itemize}
\label{lemma-2}
\end{lemma}


Denote $\widetilde{B}_{n}(x)\adef b_{n}(1/x)$. As in the case of
$B_{n}$ we have the following result.

\begin{lemma} $\widetilde{B}_{n}$ do not converges to $-1$ in
$L^{2}([0,1],dx)$.
\end{lemma}

\paragraph{\it Proof:} Denoting
$e_{k}(x)=\frho{\frac{1}{kx}}-\frac{1}{k}\,\frho{\frac{1}{x}}$,
we have $\beta_{p,q}(1/x)=e_{p}(x)-\frac{q}{p}\,e_{q}(x)$. Therefore,
each $\widetilde{B}_{n}$ is a (finite) linear combination of
$e_{k}$ and the statement of the lemma follows from Proposition 4.7 in
\cite{bae}.\qed


\section{Relation Between $\widetilde{B}_{n}$ and $B_{n}$}
\label{relation}

The next result is relevant in order to establish a relation between
the sequence $\{k_{i}\}$ and the square-free numbers and, on the other
hand, between $\widetilde{B}_{n}$ and $B_{n}$.

\begin{lemma}
The following conditions are equivalent

\begin{itemize}

\item[\bf{a.}]
$\displaystyle{
\sum_{j=1}^{i}\frac{\mu(k_{j})}{k_{j}}\;=\;\frac{1+b_{i}(k_{i})}{k_{i}}
}$, for $i\geqslant 2$. 

\item[\bf{b.}]
$\displaystyle{
b_{i}(x)=\sum_{j=1}^{i-1}\mu(k_{j})\,\frho{\frac{x}{k_{j}}} \; 
- \;
k_{i}\left(\sum_{j=1}^{i-1}\frac{\mu(k_{j})}{k_{j}}\right)\frho{\frac{x}{k_{i}}}
}$,
for $i\geqslant 2$. 

\item[\bf{c.}]
$\displaystyle{
\frac{\mu(k_{i})}{k_{i}}=\frac{1+b_{i}(k_{i})}{k_{i}}\;
- \; \frac{1+b_{i-1}(k_{i-1})}{k_{i-1}}
}$, for $i\geqslant 3$. 

\end{itemize}

\noindent
Furthermore, if the condition $k_{i+1}<2k_{i}$ is valid for $i\geqslant 2$, then all conditions
above are also equivalent to the following ones

\begin{itemize}

\item[\bf{d.}]
$\displaystyle{
\mu(k_{i+1})=b_{i}(k_{i+1}) \; - \; b_{i}(k_{i})
}$, for $i\geqslant 2$. 

\item[\bf{e.}]
$\displaystyle{
\sum_{j=1}^{i}\mu(k_{j})\pent{\frac{k_{i}}{k_{j}}}=1
}$, for $i\geqslant 1$. 

\end{itemize}
\label{lemma-3}
\end{lemma}

Comparing (\ref{eq-7}) and expression in Lemma \ref{lemma-3} (b) we can
state the following conjecture.

\begin{conjecture}
$\widetilde{B}_{n}\stackrel{?}{=}B_{n}$, for all $n\in\nn$. 
\label{conj-1}
\end{conjecture}

Observe that $\widetilde{B}_{n}$ is not a subsequence of $B_{n}$,
strictly speaking. Sequence $\{\widetilde{B}_{i}\}$ depends on
the numbers $\kk\adef\{k_{i}\}_{i\in\nn}$, the firsts of them are
given by
$$
\kk=\{1,2,3,5,6,7,10,11,13,14,15,\dots\}.
$$
These are all square-free numbers (incidentally, we prefer to
denominate the numbers in $\mm\adef\{k\in\nn\,:\,\mu(k)\neq 0\}$ as
\emph{M\"obius numbers} rather than ``square-free'', because they are
also cube-free, 4-th-power-free, etc.), 
and numerical evidence suggest $\kk\subseteq\mm$. 
Moreover, apparently none M\"obius
number is omitted, fact that seems to support the following 
conjecture.

\begin{conjecture}
$\kk\stackrel{?}{=}\mm$ 
\label{conj-2}
\end{conjecture}

If Conjecture \ref{conj-2} is true, then relation in Lemma
\ref{lemma-3} (e) is a well known result; see \cite[p. 66]{apo}. 
It is also known that square-free numbers are distributed with
density $6/\pi^{2}$; see \cite[Thm. 333, p. 269]{har}. 
We can estimate $k_{i+1}\approx k_{i}+\pi^{2}/6$, 
or $k_{i+1}/k_{i}\approx 1+\pi^{2}/6k_{i}$ and this
is less than $2$ for $k_{i}>\pi^{2}/6\approx 1.64$. 
Thus, condition $k_{i+1}<2k_{i}$ seems to be reasonable also. 
This highly speculative argument seems to suggest that Conjecture
\ref{conj-1} follows from \ref{conj-2}.


\section{Further Comments and Questions}
\label{conclusion}

Comparison of (\ref{eq-9}) and (\ref{eq-2}) (assuming (\ref{eq-3}))
gives 
\begin{equation}
\sum_{k=m}^{\infty}F_{n}(1/k)\left[\frac{1}{k^{s}}-\frac{1}{(k+1)^{s}}\right]
=\left(1-\zeta(s) \, \sum_{k=1}^{n}a_{k}\,\theta_{k}^{s}\right) \; 
- \; \frac{1}{m^{s}}.
\label{eq-30}
\end{equation}
For the particular case of $B_{n}$, where $m=n$, this expression is
given by  
\begin{equation}
\sum_{k=n}^{\infty}B_{n}(1/
k)\left[\frac{1}{k^{s}}-\frac{1}{(k+1)^{s}}\right]
=(s-1)\zeta(s)\int_{n}^{\infty}\frac{g(t)}{t^{s}}\,dt \; 
- \; \frac{1}{n^{s}}.
\label{eq-31}
\end{equation}
In particular,
$\modulo{g(t)}=\mathcal{O}(t^{-\frac{1}{2}})$ is a sufficient
condition for RH. An old result by de la Vall\'ee-Poussin states   
$\modulo{g(t)}=\mathcal{O}(\frac{1}{\ln t})$; see \cite[p. 92]{edw}.\\

If we apply Schwarz inequality in $l^{2}(\nn)$ to left hand
side in (\ref{eq-30}) trying to get an analog of BC for step Beurling
functions we have 
\begin{equation}
\modulo{\left(1-\zeta(s)\,\sum_{k=1}^{n}a_{k}\,\theta_{k}^{s}\right)}
\leqslant \modulo{\frac{1}{m^{s}}} \;
+\;\left(\sum_{k=m}^{\infty}\modulo{F_{n}(1/k)}^{2}\right)^{\frac{1}{2}}\,
\left(\sum_{k=m}^{\infty}\modulo{\frac{1}{k^{s}}-\frac{1}{(k+1)^{s}}}^{2}\right)^{\frac{1}{2}}
.
\label{eq-32}
\end{equation}
Now, $F_{n}(1/x)$ is a periodic function on the unbounded interval
$[1,+\infty)$, thus for any $N$ multiple of the period we have
\begin{equation}
\frac{a}{p}(N-m)\leqslant\sum_{k=m}^{N}\modulo{F_{n}(1/k)}^{2}\leqslant\frac{a}{p}N,
\end{equation}
where $p=p(n)$ is the period and
$a=a(n)\adef\sum_{k=1}^{p}\modulo{F_{n}(1/k)}^{2}$. Therefore, the series
for $\norma{F_{n}(1/k)}_{l^{2}(\nn)}$ is divergent. This argument
could explain the failure of $L^{2}$ convergence for general step
Beurling functions and it would be possible to write down along these
lines an alternative proof of Proposition 4.7 in \cite{bae}.



\vspace{1.5cm}

\noindent
Fernando Auil\\

\noindent
Instituto de F\'isica\\
Universidade de S\~ao Paulo\\  
Caixa Postal 66318\\  
CEP 05315-970\\       
S\~ao Paulo - SP\\            
Brasil\\

\noindent                  
E-mail: {\tt auil@fma.if.usp.br}


\end{document}